\documentclass[12pt]{article}

\usepackage{amsmath, amssymb, amsthm, amssymb}
\usepackage{cancel}
\usepackage{amsfonts}
\usepackage{graphicx}
\usepackage[frame,ps,matrix,arrow,curve,rotate]{xy}
\usepackage{enumerate}
\usepackage{tikz}
\usepackage{hyperref}
\usepackage{siunitx}
\usepackage{enumitem}
\usepackage{xspace}
\usepackage{centernot}
\usepackage{blkarray}
\usepackage{tikz-cd}

\usepackage{alphabeta}

\newtheorem{theorem}{Theorem}[section]
\newtheorem{lemma}[theorem]{Lemma}

\newtheorem{corollary}[theorem]{Corollary}
\newtheorem{proposition}[theorem]{Proposition}

\theoremstyle{definition}
\newtheorem{definition}[theorem]{Definition}

\theoremstyle{remark}

\newcommand{\mc}[1]{\mathcal{#1}}

\newcommand{\abs}[1]{\lvert#1\rvert}
\newcommand{\innerprod}[1]{\langle#1\rangle}

\newcommand{\set}[1]{\{#1\}}

\renewcommand{\geq}{\geqslant}

\renewcommand{\leq}{\leqslant}
\renewcommand{\emptyset}{\varnothing}

\newcommand{\liff}{\leftrightarrow}
\newcommand{\limplies}{\rightarrow}

\newcommand{\es}{\emptyset}

\DeclareMathOperator{\ord}{ord}
\DeclareMathOperator{\rad}{rad}

\newcommand{\bR}{\mathbb{R}}
\newcommand{\bC}{\mathbb{C}}

\newcommand{\bN}{\mathbb{N}}

\newcommand{\ZF}{\mathsf{ZF}}

\newcommand{\ZFC}{\mathsf{ZFC}}

\newcommand{\DC}{\mathsf{DC}}

\newcommand{\ACF}{\mathsf{ACF}}

\newcommand{\RCA}{\mathsf{RCA}}
\newcommand{\ACA}{\mathsf{ACA}}
\newcommand{\WKL}{\mathsf{WKL}}
\newcommand{\PRA}{\mathsf{PRA}}
\newcommand{\HB}{\mathsf{HB}}
\newcommand{\HN}{\mathsf{HN}}

\newcommand{\WO}{\mathsf{WO}}
\newcommand{\sfZ}{\mathsf{Z}}

\newcommand{\proves}{\vdash}
\newcommand{\nproves}{\mathbin{\cancel{\vdash}}}

\let\ol\overline

\usepackage{blindtext}
\usepackage{titlesec}

\begin{document}
    \title{A Reverse Mathematical Analysis of Hilbert's Nullstellensatz and Basis Theorem}
    \author{Dhruv Kulshreshtha} 
    \maketitle

\begin{abstract} This paper presents an expository reverse-mathematical analysis of two fundamental theorems in commutative algebra: Hilbert's Nullstellensatz and Basis Theorem. In addition to its profound significance in commutative algebra and algebraic geometry, the Basis Theorem is also historically notable for its nonconstructive proof. The Nullstellensatz, on the other hand, is noteworthy as it establishes a fundamental connection between the more algebraic notion of ideals
and the more geometric notion of varieties.

We explore the conscious shift from computational to conceptual approaches in mathematical argumentation, contextualizing Hilbert's contributions. We formalize the relative constructivity of these theorems using the framework of reverse mathematics, although we do not presuppose familiarity with reverse mathematics. Drawing from contemporary mathematical literature, we analyze the Basis Theorem's reliance on nonconstructive methods versus the more constructive nature of the Nullstellensatz.

Our study employs the standard tools of reverse mathematics, in particular subsystems of second-order arithmetic, to outline the minimal set-existence axioms required for these theorems. We review results showing that certain formulations of the Nullstellensatz are provable in the weak axiom system of $\RCA_0$, while the Basis Theorem requires stronger axioms, such as $\Sigma^0_2$-Induction. Consequently, we position these theorems separately within the Friedman-Simpson hierarchy. This analysis contributes to a deeper understanding of the foundational requirements for these pivotal results in algebra.
\end{abstract}

\section{Introduction}
The Hilbert Basis Theorem, proven by David Hilbert in his seminal 1890 paper on invariant theory \cite{Hilbert1890}, states that a polynomial ring in one variable -- and thus, by induction, in finitely many variables -- over a Noetherian\footnote{Note that the Basis Theorem was formulated very differently in Hilbert's original paper. In particular, the now commonly-used terminology of ``Noetherian" rings was first introduced by Chevalley in 1943 \cite{Chevalley1943}, after Noether's landmark 1921 paper \cite{Noether1921}.} ring (see Definition \ref{def: Noetherian Ring}) is itself Noetherian. In addition to its mathematical importance in commutative algebra and algebraic geometry, Hilbert's nonconstructive\footnote{We mean ``nonconstructive" in the minimal sense, i.e. proving existence without the means to construct an example.} proof is also of historical significance, representing a broader conscious\footnote{See, for example, Hilbert's quote in \cite[Preface, p.X]{Hilbert1998}: ``I have tried to avoid Kummer’s elaborate computational machinery so that here
too Riemann’s principle may be realized and the proofs compelled not by calculations but by
thought alone."} shift in mathematical thinking from computation to conception.

Before the more famous, long-lasting, and polarizing debate on Zermelo's Axiom of Choice, which Hilbert described as the axiom ``most attacked up to the present in the mathematical literature" (\cite{Hilbert1926}; see \cite[Introduction]{Moo82}), the Basis Theorem had its very own controversy. Upon learning that Hilbert could prove that any ideal in the polynomial ring was finitely generated, yet could not explicitly write down its generators, Paul Gordan, then at the forefront of invariant theory, was reliably reported to have stated ``Das ist nicht Mathematik. Das ist Theologie," i.e. ``This is not mathematics. This is theology." In fairness to Gordan, it is worth noting that he later became a major proponent and developer of Hilbert’s ideas.\footnote{As with many thumbnail stories in the history of mathematics, the real situation is more complex and interesting: McLarty accounts for numerous perspectives on the interpretation of this remark, and provides a mathematical description of the nonconstructive content in the Basis Theorem, in \cite{Mclarty2012}.}

Hilbert, however, was neither alone nor the first on the ``conceptual" side of this debate: Riemann\footnote{A key text for the development of the Riemannian style was his profound PhD thesis on the theory of complex functions \cite{Riemann1851}. A brief overview of the thesis, setting it in context and charting its significance, is in \cite[Chapters 1.2 and 1.3]{Laugwitz1999}.}, in 1850s G\"ottingen, inaugurated this broader debate against the purely ``computational" approach of mathematical argumentation associated with Kummer, Kronecker and Weierstrass, among other Berlin-mathematicians.\footnote{Notably, such conceptual argumentation can be traced even further back to G\"ottingen-mathematicians such as Dirichlet and Gauss (as highlighted in \cite{AM2014} and \cite{Ferreiros2007}).} The turn of the century also saw several other mathematicians, such as Cantor and Dedekind -- accompanied and supported\footnote{Consider Hilbert's description of Cantor's new transfinite arithmetic: ``the most astonishing product of mathematical thought, one of the most beautiful realizations of human activity in the domain of the purely intelligible" (see \cite[XXV]{Boyer1968} for details); and his description of \cite{Dedekind1888}: ``This essay...is the most important first and profound attempt to ground elementary number theory" (see \cite[In.1 (Introduction)]{Sieg2013} for details).} by Hilbert -- drawing inspiration from the Riemannian style of thinking, and, as a result, facing opposition from those on the computational side.\footnote{A particularly interesting rivalry is that between Kronecker and Dedekind, spanning many decades: See for example their attempts at extending Kummer's theory of ideal divisors in \cite{Dedekind1879} and \cite{Kronecker1882}; and the opposing methodologies in their works on the theory of numbers in \cite{Dedekind1888} and \cite{Kronecker1887}. The former is addressed in much detail in \cite{Avigad2006} and the latter is summarized in \cite{Sieg2013}. Cantor also faced much opposition from Kronecker on the matter of transfinite numbers -- see \cite{Carey2005} for details.\label{footnote: conceptual computational}} Conceptual argumentation of this sort was even taken forward by Hilbert's collaborators and successors, including Emmy Noether -- who was trained as a computational invariant theorist under Gordan, but eventually shifted to Hilbert-style methods of mathematical argumentation.

These remarks on the history of this computation/conception debate are relevant, at least spiritually, to this paper; and their relevance will hopefully become more clear to the reader once we begin formally analyzing how nonconstructive certain standard mathematical results are. What it means for a result to be more/less nonconstructive than another is addressed later in the introduction and formalized in Section \ref{Sec: Bg Rev Math}. Furthermore, an explicit example involving the relative nonconstructive content in closely related results about algebraic closures will be given in Section \ref{sec: algebraic closure}, and analogous results in the case of real closures will be described in Section \ref{sec: real closures}.

For now, we resume our historical discussion with the goal of introducing the second result in the title of this paper. Since a more thorough overview of the computation/conception debate is not directly within the scope of this paper, interested readers may refer to the following more detailed sources: Edwards provides historical context on the broader nineteenth-century shift toward conceptual argumentation in \cite{Edwards1980}; Beeson also summarizes this shift in \cite[Historical Appendix]{Beeson1985}; Rowe provides a brief historical overview of the Berlin-G\"ottingen rivalry in \cite{Rowe2000}; Avigad, as mentioned in footnote \ref{footnote: conceptual computational}, highlights Dedekind's emphasis on conceptual over algorithmic reasoning, in his development of the theory of ideals, opposing Kronecker’s theory, in \cite{Avigad2006}; and Tappenden provides an overview of Riemann’s reorientation towards conceptual reasoning in \cite{Tappenden2006}, and the contrast between the Riemannian and Weierstrassian approaches to elliptic functions and their generalisations in \cite{Tappenden2023}. Finally, this debate is also addressed in \cite{Sieg2013}, excellently summarized as: ``At the heart of the difference between these foundational positions is the freedom of introducing abstract concepts -- given by \textit{structural definitions}" \cite[In.1 (Introduction)]{Sieg2013}.

In 1893, Hilbert published his second paper on invariant theory \cite{Hilbert1893}, in which he proved the Nullstellensatz\footnote{The now commonly-used terminology of ``Nullstellensatz" was first introduced by van der Waerden in 1926 \cite{VanDerWaerden1926}.}, or zero-locus-theorem, which established a fundamental connection between the more algebraic notion of \textit{ideals} and the more geometric notion of \textit{varieties}. More specifically, ideals (see Definition \ref{def: ideal}) are more algebraic in the sense that they can be seen as generalizing the properties of prime numbers; and varieties (see Definition \ref{def: variety (common zero set)}) are more geometric in the sense that ``they are sets of points in space, generalizing geometric objects like circles, lines, and spheres, that are described by polynomial equations," as described by Karen Smith in an enlightening conversation.

One formulation of the Nullstellensatz goes as follows: Fix a field $F$, an algebraically closed field extension $K$ of $F$, and an ideal $I \subseteq F[x_1,\ldots,x_n]$. If $p \in F[x_1,\ldots,x_n]$ vanishes on $\mathbb{V}(I)$ then $p^r \in I$ for some $r\in \bN_{>0}$, where $\mathbb{V}(I)$ is the variety determined by $I$, defined in Definition \ref{def: variety (common zero set)}. An immediate corollary of the Nullstellensatz is the so-called Weak Nullstellensatz which states that if $I \subseteq F[x_1,\ldots,x_n]$ is a proper ideal, then $\mathbb{V}(I) \neq \es$, i.e. there exists a common zero for all the polynomials in the ideal in every algebraically closed extension of $F$. This may explain the name of the Nullstellensatz -- which can be proven from the Weak Nullstellensatz using the Rabinowitsch trick \cite{Rabinowitsch1930}: a simple, but clever algebraic trick that translates the weak version to the full version by a series of algebraic manipulations, following the introduction of an auxiliary variable. It is worth explicitly noting that since the Rabinowitsch trick boils down to simple algebraic manipulations, this transition is constructive.

Hilbert's first proof of the Nullstellensatz uses the Basis Theorem; however, the  Nullstellensatz is seemingly more constructive, at least in its original formulation, than the Basis Theorem. Many authors have even explicitly highlighted the constructive nature of the Nullstellensatz -- see, for example, \cite{Arrondo2006} and \cite{Manetti2005} -- and the nonconstructive nature of the Basis Theorem -- see, for example, \cite{Simpson1988}, \cite{Mclarty2012} and \cite{Stevens2023}. We formally describe this apparent difference in constructivity using the tools of \textit{reverse mathematics}.

In Section \ref{Sec: Bg Rev Math}, we introduce the research program of reverse mathematics, introduced by Harvey Friedman (see \cite{Friedman1975}) -- and developed by several others -- with the goal of determining the minimal set-existence axioms that are needed to prove theorems of non-foundational mathematics (see Section \ref{sec: subsystems of Z2}). Using reverse mathematical techniques, one may argue about theorems of non-foundational mathematics from a more foundational perspective.\footnote{As evidenced by the use of reverse mathematics in studying the logical strength of theorems in topological dynamics, combinatorics, measure theory, commutative algebra, and so on -- see, for example \cite{BlassHirstSimpson1987}, \cite{BrownSimpson1986}, \cite{BrownSimpson1993},\cite{FriedmanHirst1991}, \cite{FSS1983}, \cite{Hirst1987}, \cite{HatzikiriakouSimpson1989}, \cite{Simpson2015}, \cite{Simpson1988(B)}, and \cite{YuSimpson1990}.} For example, one may prove that some such result is independent\footnote{The logic community is not unfamiliar with independence proofs involving set-theoretic statements, like the famous results of Kurt G\"odel and Paul Cohen (\cite{Godel1938} and \cite{Cohen1963}; see \cite{Kunen1980} for accessible proofs). Reverse mathematics, however, provides the machinery to argue the independence of results in, for example, measure theory or commutative algebra.} of certain systems of mathematics; and one may even ask whether the minimal axioms that are necessary to prove the Heine-Borel Covering Lemma are also necessary/sufficient to prove the Bolzano-Weierstrass Theorem (see Theorem \ref{thm: results not provable in RCA_0}), and vice versa. 

In particular, in Section \ref{Sec: Bg Rev Math}, we introduce subsystems of second-order arithmetic\footnote{For now, we will refer to some of these systems without explaining them. Note, however, that we are not presupposing any background knowledge in reverse mathematics from our readers, and will be defining these systems in Section \ref{Sec: Bg Rev Math}.}, in increasing order of their relative strengths -- which we can then use to formally argue the relative nonconstructive content in the aforementioned kinds of results. However, there are some non-trivial challenges in such an analysis of Hilbert's Nullstellensatz and Basis Theorem.

For example, in the reverse mathematical context, and more broadly, there have been numerous formulations of the Nullstellensatz.\footnote{See, for example, \cite[Chapter 7]{AtiyahMacDonald1969}, \cite[Theorem 4.19]{Eisenbud1995}, \cite[\S 2]{FSS1983}, \cite[Chapter 1]{JL1989}, \cite[Theorem 1]{GallegoCGomezRamirez2019}, \cite[\S 1]{SakamotoTanaka2004}, \cite[\S 2.2]{Shafarevich}, and \cite[\S 2.3]{KarenSmith}.\label{footnote: formulations of Nullstellensatz}} In $\mathsf{ZFC}$, these are all provable, hence provably equivalent. However, in weak subsystems of second-order arithmetic, equivalences, or even implications, among these formulations are not so clear.\footnote{In an insightful email conversation, John Baldwin mentioned that a standard formulation of Hilbert's Nullstellensatz is equivalent to the model completeness of $\mathsf{ACF}$ (the theory of algebraically closed fields), so it should be low in the (Friedman-Simpson) hierarchy (see Section \ref{Sec: Bg Rev Math}). Furthermore, Matthew Harrison-Trainor mentioned that there could be a difference of where Hilbert's Nullstellensatz would fit in the hierarchy depending on on whether it is coded using ideals, generators, or just subsets of the field. This remark highlights an important fact to remember: Reverse
mathematical results may depend on the manner in which statements are formulated set-theoretically.}

In Sections \ref{sec: Nullstellensatz} and \ref{sec: Basis Theorem}, we describe some formulations of the Nullstellensatz and Basis Theorem, as coded in the reverse mathematics literature, along with where they fit in the Friedman-Simpson hierarchy. More specifically, in \cite{FSS1983} and \cite{SakamotoTanaka2004}, Friedman, et al.\ and Sakamoto and Tanaka show that certain formulations of the Nullstellensatz are provable in the subsystem $\RCA_0$ of second-order arithmetic (see Section \ref{sec: Nullstellensatz}). On the other hand, in \cite{Simpson1988}, Simpson shows that the Basis Theorem is not provable in $\RCA_0$. Furthermore, it is a consequence of \cite[Theorem 2.2]{Simpson2015} that the Basis Theorem is provable in $\RCA_0 + \Sigma^0_2$-Induction (see Section \ref{Sec: Bg Rev Math}), which provides an upper bound for its nonconstructive content -- which, as we will see, comes from what one would expect to be a straightforward inductive argument.

\textit{Acknowledgements.} The content of this paper is based largely on the author's undergraduate honors thesis, advised by Jamie Tappenden and Andreas Blass, whose unwavering support and invaluable feedback he is eternally grateful for. The author also thanks Jeremy Avigad, John Baldwin, Hanna Bennett, Ronnie Chen, Stephen DeBacker, Matthew Harrison-Trainor, Paul Kessenich, Sarah Koch, Jeff Lagarias, Dave Marker, Colin McLarty, Scott Schneider, Wilfried Sieg, Karen Smith, and Andrew Snowden for their insights, encouragement, and mentorship.

\section{Background: Commutative Algebra}\label{Sec: Bg Comm Alg}
We will presuppose familiarity with the basic definitions of abstract algebra (group, ring, field, etc.) as presented in standard textbooks such as \cite{AtiyahMacDonald1969} and \cite{Eisenbud1995}. Henceforth, we use the term \textit{ring} to mean a unital commutative ring, i.e. a ring with a multiplicative identity in which multiplication commutes.

\begin{definition}[Ideal]\label{def: ideal} Given a ring $R$, an \textit{ideal} $I$ of $R$ is a nonempty subset of $R$ that satisfies the following
\begin{enumerate}
    \item If $x \in I$ and $y\in I$, then $x-y \in I$. 
    \item If $x \in I$ and $a \in R$, then $ax \in I$.
\end{enumerate}
\end{definition}

\begin{definition}\label{def: proper, prime, maximal ideal} Fix a ring $R$ and an ideal $I \subseteq R$. $I$ is said to be
\begin{enumerate}
    \item \textit{Proper}, if $I \subsetneq R$.
\end{enumerate}
If, furthermore, $I$ is a proper ideal of $R$, then $I$ is said to be
\begin{enumerate}\setcounter{enumi}{1}
    \item \textit{Prime}, if for any $x,y \in R$, if $xy \in I$ then $x \in I$ or $y \in I$.
    \item \textit{Maximal}, if there is no proper ideal $J$ of $R$ such that $I \subsetneq J$.
\end{enumerate}
\end{definition}

\begin{definition}[Noetherian Ring]\label{def: Noetherian Ring} A ring $R$ is said to be \textit{Noetherian} if (any of) the following equivalent\footnote{These conditions are certainly equivalent over $\ZF+\DC$, where $\ZF$ is Zermelo-Fraenkel set theory (see any standard set theory textbook; for example, \cite[Chapter 2]{Enderton1977}), and $\DC$ is the principle of dependent choice (see, for example, \cite[Form 43]{HR98}).} conditions hold:
\begin{enumerate}
    \item Every ideal in $R$ is finitely generated.
    \item Every non-empty set of ideals in $R$ has a maximal element.
    \item The ascending chain condition holds for ideals in $R$. That is, given ideals $I_n$ in $R$ such that for each $n$, $I_n \subseteq I_{n+1}$, then there is $m$ such that for each $j\geq m$, $I_j = I_{j+1}$.
\end{enumerate}
\end{definition}

\begin{theorem}[Basis Theorem \cite{Hilbert1890}]\label{thm: Basis Theorem og} If $R$ is Noetherian, then so is $R[x]$, i.e. the one-variable polynomial ring over $R$.
\end{theorem}

A standard proof of the Basis Theorem can be found in \cite[Theorem 1.2]{Eisenbud1995} and \cite[Theorem 7.5]{AtiyahMacDonald1969}. The following consequence is the result of a straightforward (over $\ZFC$, for example) inductive argument.\footnote{We will see in Section \ref{sec: Basis Theorem} that over weak subsystems of second-order arithmetic, this argument is not quite as straightforward as one would hope.}

\begin{corollary} If $R$ is Noetherian, then so is $R[x_1,\ldots,x_n]$, i.e. the $n$-variable polynomial ring over $R$.
\end{corollary}

Fix a field $F$, and an algebraically closed\footnote{Recall that field $K$ is said to be \textit{algebraically closed} if for any non-constant polynomial $h(x) \in K[x]$, there exists $c \in K$ such that $h(c)=0$.} field extension $K$ of $F$.

\begin{definition}\label{def: variety (common zero set)} Given $S \subseteq F[x_1,\ldots,x_n]$, define the \textit{variety} determined by $S$, denoted $\mathbb{V}(S)$, to be the set of all common zeros of $S$, i.e. \[\mathbb{V}(S) = \set{\ol{x} \in K^n \mid \forall f \in S\ (f(\ol{x}) = 0)}.\]
\end{definition}

Then Hilbert's Nullstellensatz can be formulated as follows.
\begin{theorem}[Nullstellensatz \cite{Hilbert1893}]\label{Thm: Nullstellensatz og} Let $F$ and $K$ be as above, and $I \subseteq F[x_1,\ldots,x_n]$ be an ideal. If $p \in F[x_1,\ldots,x_n]$ vanishes on $\mathbb{V}(I)$, i.e. $p(\ol{x}) = 0$ for every $\ol{x} \in \mathbb{V}(I)$, then $p^r \in I$ for some $r \in \bN_{>0}$.\footnote{It is worth explicitly noting that if $I = \innerprod{g_1,\ldots,g_k}$, i.e. the ideal generated by (smallest ideal containing) $g_1,\ldots,g_k$, then $\mathbb{V}(I) = \mathbb{V}(\set{g_1,\ldots,g_k})$.}
\end{theorem}

A standard algebraic proof of the Nullstellensatz (formulated slightly differently) can be found in \cite[Theorem 4.19]{Eisenbud1995}; and a standard model-theoretic proof can be found in \cite[Theorem 3.2.11]{Marker2010}.

An immediate corollary is the so-called Weak Nullstellensatz, which can be formulated as follows.

\begin{corollary}[Weak Nullstellensatz]\label{Thm: Nullstellensatz weak og} The ideal $I \subseteq F[x_1,\ldots,x_n]$ contains $1$ if and only if the polynomials in $I$ have no common zeros in $K^n$.
\end{corollary}

An equivalent formulation states that if $I \subseteq F[x_1,\ldots,x_n]$ is a proper ideal, then $\mathbb{V}(I) \neq \es$, i.e. there exists a common zero for all the polynomials in the ideal in every algebraically closed extension of $F$. As stated previously, the Nullstellensatz can be proven from the Weak Nullstellensatz using the Rabinowitsch trick \cite{Rabinowitsch1930}. As remarked earlier, the Rabinowitsch trick boils down to simple algebraic manipulations; hence, this transition does not increase the non-constructive content.

\section{Background: Reverse Mathematics}\label{Sec: Bg Rev Math}
\subsection{Second Order Arithmetic}
The axiomatic system $\sfZ_2$ of Second Order Arithmetic has the language $\mc{L}_2$, defined below.

\begin{definition}[$\mc{L}_2$]\label{def: language of Z2} The language $\mc{L}_2$ of second-order arithmetic is a two-sorted language. This
means that there are two distinct sorts of variables which are intended
to range over two different kinds of objects. Variables of the first sort
are known as \textit{number variables}, are denoted by $k,m, n,\ldots$, and are
intended to range over the set $\omega = \set{0,1,2,\ldots}$ of all natural numbers.\footnote{Here, we use $\omega$, rather than $\mathbb{N}$, to denote the set $\set{0,1,2,\ldots}$ of natural numbers because that is the more commonly used notation in logical/set-theoretic contexts. In some cases, when we present algebraic results (such as Theorem \ref{Thm: Nullstellensatz og}) as they occur ``in nature," so to speak, we use $\mathbb{N}$ to denote the set of natural numbers and $\mathbb{N}_{>0}$ to denote the set $\bN\setminus\set{0}$. Finally, on some occasions (such as Theorem \ref{thm: provable ordinals} and Lemma \ref{lem: RCA0 not prove WO omega^omega}), $\omega$ is used (ambiguously) not to signify the set of natural numbers but rather the order type of the set of natural numbers ordered by the usual $<$ relation.} Variables of the second sort are known as \textit{set variables}, are denoted by $X,Y,Z,\ldots$, and are intended to range over all subsets of $\omega$. $\mc{L}_2$ has countably many set and number variables.

\textit{Numerical terms}, intended to denote natural numbers, are number variables, the constant symbols
$0$ and $1$, and $t_1 + t_2$ and $t_1 \cdot t_2$\footnote{Strictly speaking, these should be $(t_1+t_2)$ and $(t_1\cdot t_2)$ to avoid ambiguity in the order of operations, but we take such details for granted.}, whenever $t_1$ and $t_2$ are numerical terms. Here $+$ and $\cdot$ are binary operation symbols intended to denote addition and multiplication of natural numbers.

\textit{Atomic formulas} are $t_1 = t_2$, $t_1 < t_2$, and $t_1 \in X$ where $t_1$ and $t_2$ are numerical terms and $X$ is any set variable, and are intended to mean that $t_1$ equals
$t_2$, $t_1$ is less than $t_2$, and $t_1$ is an element of $X$ respectively.

\textit{Formulas} are built up
from atomic formulas by means of propositional connectives $\land, \lor, \neg, \limplies$ (and, or, not, if...then), and quantifiers $\forall$ and $\exists$ (representing the universal and existential quantifiers). There are, however, two sorts of each quantifier: \textit{numerical quantifiers} (denoted $\forall n, \exists n$ with lower case variables) and \textit{set quantifiers} (denoted $\forall X, \exists X$ with upper case variables). 
\end{definition}

The language $\mc{L}_1$ for First Order Arithmetic is just $\mc{L}_2$ without set quantifiers or set variables. An $\mc{L}_2$-formula is said to be \textit{arithmetical} if it contains no set quantifiers, i.e. all the quantifiers appearing are number quantifiers.

\begin{definition}[$\sfZ_2$]\label{def: Second Order Arithmetic}
The axioms of Second Order Arithmetic ($\sfZ_2)$ are as follows:
\begin{enumerate}
    \item Basic Axioms: For all natural numbers $n,m$,
    \begin{itemize}
        \item $m+1 \neq 0$
        \item $m+1=n+1 \limplies m=n$
        \item $m+0 = m$
        \item $m+(n+1) = (m+n)+1$
        \item $m\cdot 0 = 0$
        \item $m\cdot(n+1) = (m\cdot n)+m$
        \item $\neg(m<0)$
        \item $m<n+1 \liff (m<n) \lor (m=n)$
    \end{itemize}
    \item Induction Axiom:\hspace{-1mm} $\forall X[(0 \in X \land \forall n(n\in X\hspace{-0.5mm} \limplies\hspace{-0.5mm} n+1\in X)) \hspace{-0.5mm}\limplies\hspace{-0.5mm} \forall n(n\in X)]$
    \item Comprehension Scheme: $\exists X \forall n (n\in X \liff \varphi(n))$ where $\varphi(n)$ is any $\mc{L}_2$-formula in which $X$ doesn't occur freely.
\end{enumerate}
It is a consequence of (2) and (3) that $\sfZ_2$ implies the full second order induction scheme: For any $\mc{L}_2$-formula $\varphi,$
\[ (\varphi(0) \land \forall n(\varphi(n)\limplies\varphi(n+1)))\limplies\forall n\varphi(n).\]
\end{definition}

\begin{definition}[Models of $\sfZ_2$] 
An $\mc{L}_2$-structure is an ordered 7-tuple
\[ M = \innerprod{\abs{M},\mc{S}_M, +_M, \cdot_M, <_M, 0_M, 1_M},\]
where $\abs{M}$ is a set which serves as the domain of the number variables, $\mc{S}_M$ is a set of subsets of $\abs{M}$ serving as the domain of the set variables, $+_M$ and $\cdot_M$ are binary operations on $\abs{M}$, $0_M$ and $1_M$ are distinguished elements of
$\abs{M}$, and $<_M$ is a binary relation on $\abs{M}$. We always assume that the sets $\abs{M}$ and $\mc{S}_M$ are disjoint and nonempty. Formulas of $\mc{L}_2$ are interpreted in $M$ in the obvious way. A \textit{model} of $\sfZ_2$ is an $\mc{L}_2$-structure that satisfies all the axioms in Definition $\ref{def: Second Order Arithmetic}$.

For the purposes of this paper, the domain of numbers $\abs{M}$, and the operations and relations $(+_M,\cdot_M,<_M)$, can be thought of as the usual operations of natural numbers (see $\omega$-models in \cite[Chapter VIII]{Simpson2009}).
\end{definition}

\subsection{The Arithmetical and Analytical Hierarchies}\label{sec: arithmetical hierarchy}

Given a formula $\varphi$, a number variable\footnote{The notions of bound and free can all be defined analogously for set variables.} $n$ is said to be \textit{bound} in $\varphi$ if every occurrence of $n$ in $\varphi$ is within the scope of a quantifier of $n$; and $n$ is said to be \textit{free} if it is not bound.\footnote{A formula with no free variables is called a \textit{sentence}.} For example, in $\varphi(n) := \forall n(n+n=0)$, $n$ is bound by the quantifier $\forall n$, and in $\psi(n,m) := \exists m (n+m = 0)$, $n$ is free.

Let $t$ be a term that does not contain $n$. We abbreviate $\forall n(n<t\ \limplies \varphi)$ as $(\forall n<t)\varphi$ and $\exists n(n<t\ \land \varphi)$ as $(\exists n<t)\varphi$. The quantifiers $\forall n<t$ and $\exists n<t$ are called \textit{bounded quantifiers}. A \textit{bounded quantifier formula} is a formula whose quantifiers are all bounded number
quantifiers.

\begin{definition}[Arithmetical Hierarchy]\label{def: Arithmetical Hierarchy} For $k\in \omega = \set{0,1,2,\ldots}$, an $\mc{L}_2$-formula $\varphi$ is said to be $\Sigma^0_k$ if it is of the form
\[ \exists n_1 \forall n_2 \exists n_3 \cdots n_k\ \theta, \]
and an $\mc{L}_2$-formula $\varphi$ is said to be $\Pi^0_k$ if it is of the form
\[ \forall n_1 \exists n_2 \forall n_3 \cdots n_k\ \theta, \]
where $n_1,\ldots,n_k$ are are number variables
and $\theta$ is a bounded quantifier formula. An $\mc{L}_2$-formula is said to be \textit{arithmetical}, if it is $\Sigma^0_k$ or $\Pi^0_k$ for some $k \in \omega$.
\end{definition}

In both cases, $\varphi$ consists of $k$ alternating unbounded number quantifiers followed by a formula containing only bounded number quantifiers. In the $\Sigma^0_k$ case, the first unbounded number quantifier is existential, while in the $\Pi^0_k$ case it is universal (assuming $k \geq 1$). Thus for instance a $\Pi^0_2$ formula is of the form $\forall m \exists n\ \theta$,
where $\theta$ is a bounded quantifier formula. A $\Sigma^0_0$ or $\Pi^0_0$ formula is the same
thing as a bounded quantifier formula.

There is no need to include, in Definition \ref{def: Arithmetical Hierarchy}, additional clauses covering the cases of multiple non-alternating quantifiers, such as $\exists n_1 \forall n_2 \forall n_3 \exists n_4\ \theta$, since these will always be equivalent to formulas where all the unbounded quantifiers alternate. For example, $\exists n_1 \forall n_2 \forall n_3 \exists n_4\ \theta$ is equivalent to \[\exists n_1\ \forall m\ \forall n_2\hspace{-0.75mm}<\hspace{-0.75mm}m\ \forall n_3\hspace{-0.75mm}<\hspace{-0.75mm}m\  \exists n_4\  \theta,\]
which can easily be translated to a $\Sigma^0_3$ formula using a sequence coding function, as described in \cite[Chapter 6]{Shoenfield1967}.

Clearly any $\Sigma^0_k$ formula is logically equivalent to the negation of a $\Pi^0_k$ formula, and vice versa. Moreover, up to logical equivalence of formulas, we have $\Sigma^0_k \cup \Pi^0_k \subseteq \Sigma^0_{k+1} \cap \Pi^0_{k+1}$ for each $k \in \omega$. We say an $\mc{L}_2$-formula $\varphi$ is $\Delta^0_k$ if it is equivalent to both a $\Sigma^0_k$ formula and a $\Pi^0_k$ formula. We make the properties of this hierarchy more explicit as follows.

\begin{proposition}[See Sections I.7 and IX.1 of \cite{Simpson2009}] The following hold (up to logical equivalence of formulas):
\begin{enumerate}
    \item $\Pi^0_k \subsetneq \Pi^0_{k+1}$ and $\Sigma^0_k \subsetneq \Sigma^0_{k+1}$, for any $k \in \omega$.
    \item $\Delta^0_k \subsetneq \Sigma^0_k, \Pi^0_k$, for any $k \geq 1$.
    \item $\Sigma^0_k \cup \Pi^0_k \subsetneq \Delta^0_{k+1}$, for any $k \geq 1$.
\end{enumerate}
\end{proposition}

These results (for $k\geq 1$) can be summarized using the following figure, where the arrows indicate strict inclusions.
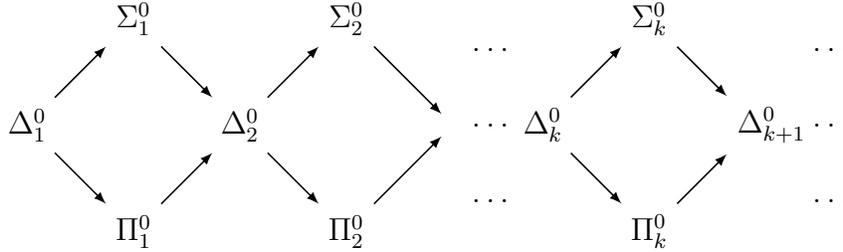
\begin{center}
\begin{tikzpicture}[node distance=2cm,semithick,-latex]
\node(1){$\Delta_1^0$};
\node(1a)[above right of=1]{$\Sigma_1^0$};
\node(1b)[below right of=1]{$\Pi_1^0$};
\node(2)[above right of=1b]{$\Delta_2^0$};
\path (1) edge (1a);
\path (1) edge (1b);
\path (1a) edge (2);
\path (1b) edge (2);

\node(2a)[above right of=2]{$\Sigma_2^0$};
\node(2b)[below right of=2]{$\Pi_2^0$};
\node(3)[above right of = 2b]{};
\path (2) edge (2a);
\path (2) edge (2b);
\path (2a) edge (3);
\path (2b) edge (3);

\node(3)[right of = 3,node distance=0.5cm]{$\ldots$};
\node [above of=3,node distance = 1cm]{$\ldots$};
\node [below of=3,node distance = 1cm]{$\ldots$};

\node(3)[right of = 3, node distance = 0.7cm]{$\Delta_k^0$};

\node(3a)[above right of=3]{$\Sigma_k^0$};
\node(3b)[below right of=3]{$\Pi_k^0$};
\node(4)[below right of=3a]{$\quad\Delta_{k+1}^0$};
\path (3) edge (3a);
\path (3) edge (3b);
\path (3a) edge (4);
\path (3b) edge (4);

\node (5) [right of=4,node distance = 1cm]{$\ldots$};
\node (5a) [above of=5,node distance = 1cm]{$\ldots$};
\node (5b) [below of=5,node distance = 1cm]{$\ldots$};
\end{tikzpicture}
\end{center}

In parallel, the analytical hierarchy, which is the extension of the arithmetical hierarchy with set quantifiers, can be defined similarly. In this case, we quantify over sets of numbers and use arithmetical formulas in place of bounded quantifier formulas. We use 1 as a superscript instead of 0 to explicitly indicate the level of the hierarchy.\footnote{One can also similarly define higher levels with quantification over sets of sets of numbers, sets of sets of sets of numbers, etc., reflected in the superscripts 2, 3, etc., but an understanding of the first two levels, with superscripts 0 and 1, is sufficient for the purposes of this paper.} So for example, $\Delta^1_0 = \Sigma^1_0 = \Pi^1_0$ indicates the class of $\mc{L}_2$-formulas with number quantifiers but no set quantifiers. Furthermore, an $\mc{L}_2$-formula is $\Sigma^1_1$ if it logically equivalent to $\exists X\hspace{0.5mm} \theta$, $\Pi^1_1$ if it logically equivalent to $\forall X\hspace{0.5mm} \theta$, $\Sigma^1_2$ if it is logically equivalent to $\exists X\hspace{0.5mm} \forall Y\hspace{0.5mm} \theta$, $\Pi^1_2$ if it is logically equivalent to $\forall X\hspace{0.5mm} \exists Y\hspace{0.5mm} \theta$, and so on, where $\theta \in \Delta^1_0$.

\subsection{Subsystems of $\sfZ_2$}\label{sec: subsystems of Z2}
By a \textit{subsystem} of $\sfZ_2$, we mean a system of arithmetic with the basic axioms, i.e. \ref{def: Second Order Arithmetic}(1), and restrictions of induction or comprehension, i.e. \ref{def: Second Order Arithmetic}(2 and 3). For this section, and the rest of this paper, we use \cite{Simpson2009} as a reference.

We describe the subsystems $\RCA_0$, $\WKL_0$, and $\ACA_0$\footnote{These subsystems are three of the ``Big Five" subsystems of the Friedman-Simpson hierarchy -- as described in \cite{NS2023}.} of $\sfZ_2$, which differ in their set existence axioms. Furthermore, $\RCA_0$ is no stronger than $\WKL_0$ in terms of provability, which again is no stronger than $\ACA_0$ in this sense (see \cite[I.10.2]{Simpson2009}). That is,
\[ \RCA_0 \vdash \varphi \implies \WKL_0 \vdash \varphi \implies \ACA_0 \vdash \varphi. \]
Soon we will see that these implications are not reversible, hence $\RCA_0$ is strictly weaker than $\WKL_0$, which again is strictly weaker than $\ACA_0$.

We first describe $\PRA$, the formal system of primitive recursive arithmetic, which is viewed by many as a plausible explication of ``finitistically acceptable reasoning" -- see, for example, \cite{Tait1981}\footnote{``We shall see that there is no question but that [primitive recursive]
reasoning is finitist" \cite{Tait1981}.} and \cite{Simpson1988(B)}. In Section \ref{sec: conservativity and hilbert's program}, we state a conservativity result, due to Harrington, for $\WKL_0$ over $\PRA$, which Simpson argues is a partial realization of \textit{Hilbert's program} for the foundations of mathematics (see \cite[Remark IX.3.18]{Simpson2009}). However, whether or not this is Hilbert's conception is disputed. For example, Sieg argues that Simpson equating finitistic reduction to Hilbert’s program is inaccurate; and that Kronecker’s name would be more appropriate to attach to said reductionist program \cite{Sieg1990}. In Section \ref{sec: conservativity and hilbert's program}, we provide some additional context on Hilbert's program, and direct interested readers to more detailed sources. 

\begin{definition}[$\PRA$]\label{def: PRA} The language of $\PRA$ is described in \cite[Definition IX.3.1]{Simpson2009}. With reference to Definition \ref{def: language of Z2}, it contains only number variables, number relations and operations from before, and symbols for each \textit{basic primitive recursive function}: the constant zero function $Z(x)=0$, the successor function $S(x)=x+1$, and the projection functions $P^k_i(x_1,\ldots,x_k)=x_i$. It also codes a symbolization of all \textit{primitive recursive functions}, which are built out of the basic functions using the composition and primitive recursion operators.

The \textit{intended model} of $\PRA$ consists of the nonnegative integers $\omega = \set{0,1,2,\ldots}$, together with the primitive recursive functions, as defined above. The axioms of $\PRA$, defined in \cite[Definition IX.3.2]{Simpson2009}, include the usual axioms for equality, zero, the successor and projection functions, composition, and so on, along with the scheme of primitive recursive induction:\medskip

\noindent \textit{Primitive Recursive Induction} \vspace{-2mm}
\[ (\theta(0)\land \forall x\ (\theta(x)\limplies \theta(\underline{S}(x))) \limplies \forall x\ \theta(x),\]
where $\theta(x)$ is any quantifier-free formula in the language of $\PRA$ with a distinguished free number variable $x$.
\end{definition}

We now describe $\RCA_0$, which is, for all intents and purposes, our base system. More specifically, equivalences such as theorems \ref{thm: WKL_0 equivalences} and \ref{thm: ACA_0 equivalences} are considered over $\RCA_0$, and algebraic definitions in, for example, Section \ref{Sec: Rev Math Analysis} are made in $\RCA_0$. 

\begin{definition}[$\RCA_0$]\label{def: RCA_0}\addcontentsline{toc}{subsubsection}{Recursive Comprehension}
Along with the basic axioms, i.e. \ref{def: Second Order Arithmetic}(1), $\RCA_0$ allows for $\Delta^0_1$-Comprehension (also called recursive comprehension) and $\Sigma^0_1$-Induction (equivalently $\Pi^0_1$-Induction), i.e.\medskip 

\noindent \textit{Recursive Comprehension} \vspace{-2mm}
\[ \forall x\ (\varphi(x) \liff \psi(x)) \limplies \exists X\ \forall x\ (x \in X \liff \varphi(x)), \]
where $\varphi$ is $\Sigma^0_1$, $\psi$ is $\Pi^0_1$, and $X$ is not free in either $\varphi$ or $\psi$. \bigskip

\noindent \textit{$\Sigma^0_1$-Induction} \vspace{-2mm}
\[ (\varphi(0) \land \forall x\ (\varphi(x) \limplies \varphi(x+1))) \limplies \forall x\ \varphi(x), \]
where $\varphi$ is $\Sigma^0_1$ (or equivalently $\Pi^0_1$).\footnote{One could analogously define, for example, $\Sigma^0_2$-Induction by allowing $\varphi$ to be $\Sigma^0_2$.\label{footnote: Sigma^0_2-induction}}
\end{definition}

The acronym $\RCA$ stands for Recursive Comprehension Axiom, and the subscript $0$ indicates that only the restricted induction schema, i.e. $\Sigma^0_1$-Induction, is being assumed. One may also consider the system $\RCA$, i.e. with the full first-order induction schema, i.e. $\Sigma^0_n$-Induction. In \cite[Introduction]{FSS1983}, Friedman, et al.\ motivate why this is not necessary, especially in the algebraic context.

Two quantifier induction, i.e. $\Sigma^0_2$-Induction or $\Pi^0_2$-Induction, is not provable in $\RCA_0$, and neither is $\Sigma^0_1$-Comprehension (see \cite[Chapter II]{Simpson2009}).

We now list some ordinary mathematical results\footnote{Carrying these over to formal arithmetic requires some coding.} that are, and aren't, provable in $\RCA_0$. We reserve results in commutative algebra for a deeper analysis in Section \ref{Sec: Rev Math Analysis}.

\begin{theorem}[See Theorem I.8.3 of \cite{Simpson2009}] The following are provable in $\RCA_0$.
\begin{enumerate}
    \item Baire Category Theorem: Let $\set{U_n : n\in \bN}$ be a sequence of dense open sets in $\bR^k$. Then there exists $x\in \bR^k$ such that $x\in U_n$ for all $n\in \bN$.
    \item Intermediate Value Theorem: If $f(x)$ is continuous on the unit interval $[0,1]$, and if $f(0)<0<f(1)$, then there exists $x \in (0,1)$ such that $f(x)=0$.
    \item Tietze Extension Theorem: Let $X$ be a complete separable metric space. Given a closed set $C\subseteq X$ and a continuous function $f: C\to [-1,1]$, there exists a continuous function $g: X\to [-1,1]$ such that $g(x)=f(x)$ for all $x\in C$.
    \item Weak version of G\"odel's Completeness Theorem: Let $X \subseteq \mathsf{Snt}$\footnote{Formally, given a countable language $\mc{L}$, i.e. a countable set of set of relation, operation, and constant symbols, we identify terms and formulas with their G\"odel numbers under a fixed G\"odel numbering, which can be constructed by
primitive recursion (see \cite[Theorem II.3.4]{Simpson2009}), using $\mc{L}$ as a parameter. We can prove in $\RCA_0$ that there exists a set $\mathsf{Snt}$ consisting of all G\"odel numbers of sentences.\label{footnote: Snt}} be consistent and closed under logical consequence. Then there exists a countable model $M$ such that $M$ satisfies $\varphi$, for all $\varphi\in X$.
    \item Soundness Theorem: If $X \subseteq \mathsf{Snt}$ and there exists a countable model $M$ such that $M$ satisfies $\varphi$, for all $\varphi\in X$, then $X$ is consistent.
\end{enumerate}
\end{theorem}

\begin{theorem}[Theorems I.9.3 and I.10.3 of \cite{Simpson2009}]\label{thm: results not provable in RCA_0} The following are not provable in $\RCA_0$
\begin{enumerate}
\item The Heine-Borel Covering Lemma: Every covering of the closed interval $[0,1]$ by a sequence of open intervals has a finite subcovering. 
\item The Bolzano-Weierstrass Theorem: Every bounded sequence of real numbers, or points in $\bR^n$, has a convergent subsequence.
\end{enumerate}
\end{theorem}

Both results in Theorem \ref{thm: results not provable in RCA_0} are non-constructive, in the sense that they guarantee the existence of certain mathematical objects, without providing an explicit algorithm to find them. In particular, \ref{thm: results not provable in RCA_0}(1) states the existence of a finite subcovering and \ref{thm: results not provable in RCA_0}(2) the existence of a convergent subsequence, without giving the machinery to determine what these actually look like or how one could compute them.

It turns out that more can be said about their non-constructive content using the framework of reverse mathematics. More specifically, \ref{thm: results not provable in RCA_0}(1) is less non-constructive, i.e. more constructive, than \ref{thm: results not provable in RCA_0}(2), since the former is provably equivalent to $\WKL_0$ (over $\RCA_0)$, and the latter is provably equivalent to the strictly stronger $\ACA_0$ (over $\RCA_0)$.

\begin{center}
\tikzset{
  point_style/.style={circle, draw, minimum size=2pt, inner sep=0pt, fill=black},
  colored_edge_1/.style={draw=black},
  colored_edge_2/.style={draw=orange},
  colored_edge_0/.style={draw=black},}
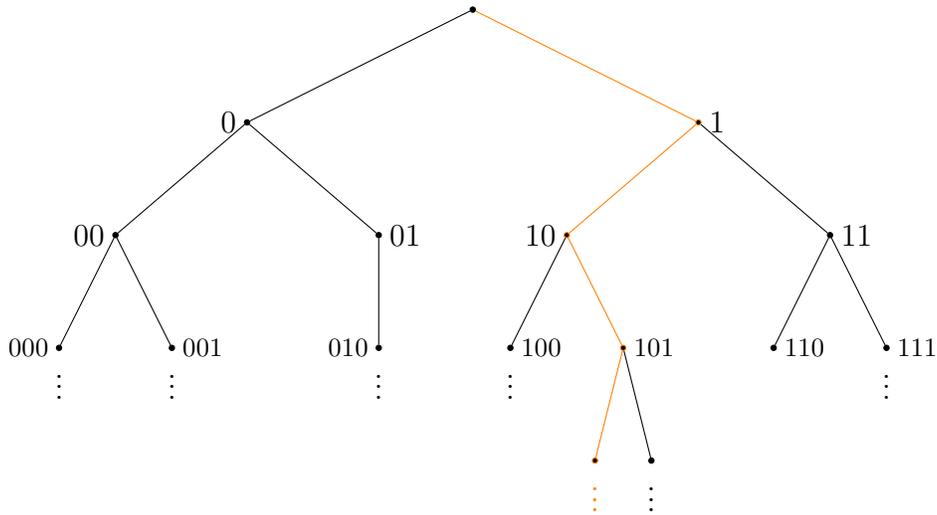
\begin{figure}[h!]
\centering
\begin{tikzpicture}[level distance=1.5cm,
  level 0/.style={sibling distance=8cm},
  level 1/.style={sibling distance=6cm},
  level 2/.style={sibling distance=3.5cm},
  level 3/.style={sibling distance=1.5cm},
  level 4/.style={sibling distance=0.75cm},
  level 5/.style={sibling distance=0.5cm},
  ]
  \node[point_style] (start) {} 
    child[colored_edge_0] {
      node[point_style] (0) {}
      child[colored_edge_0] {
        node[point_style] (00) {}
        child[colored_edge_0] {
          node[point_style] (000) {} 
          node[left] {\footnotesize 000}
          node[below] {\vdots}
        }
        child[colored_edge_0] {
          node[point_style] (001) {} 
          node[right] {\footnotesize 001}
          node[below] {\vdots}
        }
        node[left] {00}
      }
      child[colored_edge_0] {
        node[point_style] (01) {}
        child[colored_edge_0] {
          node[point_style] (010) {} 
          node[left] {\footnotesize 010}
          node[below] {\vdots}
        }
        node[right] {01}
      }
      node[left] {0}
    }
    child[colored_edge_2] {
      node[point_style] (1) {}
      child[colored_edge_2] {
        node[point_style] (10) {}
        child[colored_edge_0] {
          node[point_style] (100) {} 
          node[right] {\footnotesize 100}
          node[below] {\vdots}
        }
        child[colored_edge_2] {
          node[point_style] (101) {}
            child[colored_edge_2] {
            node[point_style] (1010) {}
            node[below] {\color{orange} \vdots}
            }
            child[colored_edge_0] {
            node[point_style] (1011) {}
            node[below] {\vdots}
            }
          node[right] {\footnotesize 101}
        }
        node[left] {10}
      }
      child[colored_edge_0] {
        node[point_style] (11) {}
        child[colored_edge_0] {
          node[point_style] (110) {} 
          node[right] {\footnotesize 110}
          node[below] {}
        }
        child[colored_edge_0] {
          node[point_style] (111) {} 
          node[right] {\footnotesize 111}
          node[below] {\vdots}
        }
        node[right] {11}
      }
      node[right] {1}
    };
\end{tikzpicture}
\caption{Binary Tree with {\color{orange} path}}
\end{figure}
\end{center}

\begin{definition}[Binary Trees]\label{def: binary trees and paths}\addcontentsline{toc}{subsubsection}{Weak K\"onig's Lemma} Let $2^{<\bN}$ (where we use the convention $2=\set{0,1}$) be the set of all finite sequences of $0$'s and $1$'s. A \textit{binary tree} is a set $T\subseteq 2^{<\bN}$ such that any initial segment of a sequence in $T$ belongs to $T$. A \textit{path}\footnote{By ``path," we always mean infinite path.} through $T$ is a function $f:\bN\to 2$ such that for every $k \in \bN$, $f[k]=(f(0),f(1),\ldots,f(k-1))$ belongs to $T$.\footnote{More generally, we may define a finitely-branching tree as a set $T\subseteq \bN^{<\bN}$ such that each node has finitely many children and any initial segment of a sequence in $T$ belongs to $T$. In this case, a path $f:\bN\to \bN$ is defined analogously.}
\end{definition}

\begin{theorem}[Weak K\"onig's Lemma (WKL)]\label{thm: WKL}
Every infinite binary tree has a path. 
\end{theorem}

See Figure 1 for an example of a path (in orange) through a binary tree, as in Theorem \ref{thm: WKL}.

\begin{definition}[$\WKL_0$]\label{def: WKL_0} $\WKL_0$ is $\RCA_0$ together with WKL.
\end{definition}

It is probably no surprise where the $\WKL$ acronym comes from. Furthermore, the subscript $0$ once again indicates that only the restricted induction schema is being assumed.

A lot of ordinary mathematics can be done in $\WKL_0$, as evidenced by the following theorem.\footnote{Due to model-theoretic extension principles such as the Upwards L\"owenheim-Skolem Theorem, the restriction to countable structures, e.g.\ countable rings, is not an actual limitation to the generality of these results. However, such machinery is outside the scope of subsystems of $\sfZ_2$.}

\begin{theorem}\label{thm: WKL_0 equivalences}[See Theorem I.10.3 of \cite{Simpson2009}] The following are provably pairwise equivalent over $\RCA_0$:
    \begin{enumerate}
        \item $\WKL_0$
        \item The Heine-Borel Covering Lemma: Every covering of the closed interval $[0,1]$ by a (countable) sequence of open intervals has a finite subcovering.
        \item Every continuous real-valued function on any compact metric space (see \cite[Definition III.2.3]{Simpson2009}) is bounded, has a supremum, and is uniformly continuous.
        \item Every continuous real-valued function on $[0,1]$ is Riemann integrable.
        \item Lindenbaum Lemma: Every countable consistent set of sentences extends to a maximal such set.
        \item G\"odel's Completeness Theorem: Every countable consistent set of sentences has a countable model.
        \item G\"odel's Compactness Theorem: Given $X \subseteq \mathsf{Snt}$ (see footnote \ref{footnote: Snt}), if each finite subset of $X$ has a model, then so does $X$.
        \item Every countable ring has a prime ideal (see Definition \ref{def: proper, prime, maximal ideal}).
        \item Brouwer’s fixed point theorem: Every continuous function $f : [0,1]^n \to [0,1]^n$ has a fixed point, i.e.\ a point $c$ such that $f(c)=c$.
    \end{enumerate}
\end{theorem}

There are still important mathematical results that $\WKL_0$ is not strong enough to prove: for example the Bolzano-Weierstrass Theorem. So we now define $\ACA_0$, which is a stronger (see \cite[I.10.2]{Simpson2009}) subsystem of $\sfZ_2$ that is strong enough to prove some of these.

\begin{definition}[$\ACA_0$]\label{def: ACA_0}\addcontentsline{toc}{subsubsection}{Arithmetical Comprehension} $\ACA_0$ is $\RCA_0$ together with arithmetical comprehension, i.e.\bigskip

\noindent \textit{Arithmetical Comprehension} \vspace{-2mm}
\[ \exists X\ \forall x\ (x \in X \liff \varphi(x)), \]
where $\varphi$ is $\Sigma^0_n$ (or equivalently $\Pi^0_n)$, for some $n$, and $X$ does not occur freely in $\varphi$.
\end{definition}

The acronym $\ACA$ stands for Arithmetical Comprehension Axiom, and the subscript $0$ serves the same purpose as before.

The following theorem motivates the relative strength of $\ACA_0$ as compared to $\WKL_0$. More specifically, since the following results are equivalent to $\ACA_0$, they are not provable in $\WKL_0$.

\begin{theorem}[See Theorem I.9.3 of \cite{Simpson2009}]\label{thm: ACA_0 equivalences} The following are provably pairwise equivalent over $\RCA_0$:
    \begin{enumerate}
        \item $\ACA_0$
        \item Bolzano-Weierstrass Theorem: Every bounded sequence of real numbers, or points in $\bR^n$, has a convergent subsequence.
        \item Every Cauchy sequence of real numbers is convergent.
        \item Monotone convergence theorem: Every bounded monotone sequence of real numbers is convergent.  
        \item Every sequence of points in a compact metric space has a convergent subsequence.
        \item Every countable ring has a maximal ideal (see Definition \ref{def: proper, prime, maximal ideal}).
        \item Every countable vector space over any countable field has a basis.
        \item Every countable field (of characteristic 0) has a transcendence basis.
        \item Every countable Abelian group has a unique (up to isomorphism) divisible closure.
        \item K\"onig's Lemma: Every infinite, finitely branching tree has a path.
    \end{enumerate}
\end{theorem}

In theorems \ref{thm: WKL_0 equivalences} and \ref{thm: ACA_0 equivalences}, we intentionally omit results pertaining to algebraic/real closures of fields. These will be covered in later sections, particularly \ref{sec: algebraic closure} and \ref{sec: real closures}.

\subsection{Conservativity results and Hilbert's program}\label{sec: conservativity and hilbert's program}
In subsequent sections, we study results involving the provability of algebraic results in certain subsystems of $\sfZ_2$, particularly $\RCA_0$ and $\ACA_0$. One may ask why authors have seemingly avoided studying some of these results in $\WKL_0$. As it turns out, $\WKL_0$ is conservative over $\RCA_0$ for the kinds of sentences that these results can be formulated in (see Theorem \ref{thm:conservativity WKL_0 over RCA_0}). Furthermore, we also study a conservativity result of $\WKL_0$ over $\PRA$\footnote{Note that although the proof of Theorem \ref{thm:conservativity WKL_0 over PRA} in \cite{Simpson2009} uses model-theoretic techniques, Sieg \cite{Sieg1985} gives a primitive recursive proof transformation which, given a proof of a $\Pi^0_2$ $\varphi$ in $\WKL_0$, generates a proof of $\varphi$ in $\PRA$. Hence, this result is provable
within a finitary system and thus allows the reduction to go through.}, and its connection (or lack thereof) to Hilbert's program.

The following theorem may be expressed by saying that $\WKL_0$ is conservative over $\RCA_0$ for $\Pi^1_1$ sentences.
\begin{theorem}[Friedman; see Corollary IX.2.6 of \cite{Simpson2009}]\label{thm:conservativity WKL_0 over RCA_0} If $\varphi$ is a $\Pi^1_1$ sentence and $\WKL_0$ proves $\varphi$, then $\RCA_0$ also proves $\varphi$.
\end{theorem}

Similarly, the following theorem may be expressed by saying that $\WKL_0$ is conservative over $\PRA$ for $\Pi^0_2$ sentences. 
\begin{theorem}[Harrington; see Theorem IX.3.16 of \cite{Simpson2009}]\label{thm:conservativity WKL_0 over PRA} If $\varphi$ is a $\Pi^0_2$ sentence and $\WKL_0$ proves $\varphi$, then $\PRA$ also proves $\varphi$.\footnote{A proof of Theorem \ref{thm:conservativity WKL_0 over PRA}, and its extension to a strengthening $\WKL_0^+$ of $\WKL_0$, is outlined in \cite[Sections 5 and 6]{Day2019}.}
\end{theorem}

In \cite{Simpson1988(B)}, Simpson argues that Theorem \ref{thm:conservativity WKL_0 over PRA} represents a partial realization of \textit{Hilbert's program} (by which he means \textit{finitistic reductionism}) for the foundations of mathematics.\footnote{As stated before, although Hilbert did not explicitly spell out a precise definition of \textit{finitism}, many  (see, for example, \cite{Simpson1988(B)} and \cite{Tait1981}) agree that the formal system of $\PRA$ (see Definition \ref{def: PRA}) captures the essence of this notion.} The goal of finitistic reductionism, in this regard, was to show that non-finitistic set-theoretical mathematics can be reduced to $\PRA$, by means of conservation results (for $\Pi^0_1$ sentences).

Of the many proof-theoretic advances of the early twentieth-century, one that is arguably the most significant\footnote{One may recall von Neumann's remark on the occasion of the presentation of the Albert Einstein Award to G\"odel in 1951: ``Kurt G\"odel’s achievement in modern logic is singular and monumental -- indeed, it is more than a monument; it is a landmark which will remain visible far in space and time."} -- as well as the most relevant to the realization (or, in this case, the lack thereof) of finitistic reductionism -- is G\"odel's publication of his incompleteness theorems \cite[Theorems VI and XI]{Godel1931}\footnote{Roughly, G\"odel's First Incompleteness Theorem \cite[Theorem VI]{Godel1931}, after incorporating subsequent refinements, is the assertion that any consistent formal system $F$ within which a certain amount of elementary arithmetic can be carried out is incomplete. That is, there are statements of the language of $F$ which can neither be proved nor disproved in $F$. Furthermore, G\"odel's Second Incompleteness Theorem \cite[Theorem XI]{Godel1931}, after incorporating subsequent refinements, is the assertion that for any consistent system $F$ within which a certain amount of elementary arithmetic can be carried out, the consistency of $F$ cannot be proved in $F$ itself.}. A consequence of these incompleteness theorems is that a complete realization of Hilbert's program is impossible.

However, it is still worth trying to understand what parts of Hilbert's program are salvageable, i.e. what sorts of infinitistic mathematics can, in fact, be reduced to finitism. Simpson (see \cite[Remark IX.3.18]{Simpson2009}) reformulates this question, in the language of subsystems of $\sfZ_2$ as: Which interesting subsystems of $\sfZ_2$ are conservative over $\PRA$ for $\Pi^0_1$ sentences? 
In this regard Theorem \ref{thm:conservativity WKL_0 over PRA} provides a partial answer to his question, and, in turn, a partial realization of Hilbert's program. This, however, relies on the assumption that finitistic reductionism is, in fact, Hilbert's conception -- which is disputed (refer to discussion in \ref{sec: subsystems of Z2}).

Since it is not directly within the scope of (but still relevant to) this paper, we pause our discussion on Hilbert's program and G\"odel's incompleteness theorems here, and refer interested readers to more detailed sources: Sieg provides a much more informative description of Hilbert's program(s) in \cite[Section II]{Sieg2013}; Eastaugh provides a much more complete overview of the partial realization story in \cite[Section 4.5]{Eastaugh2015}; Davis lists and discusses undecidability results including (and following) G\"odel's incompleteness theorems in \cite{Davis1965}; and Prince provides an easy-to-follow annotated version of G\"odel's 1931 paper in \cite{Prince2022}.

 \subsection{Provable Ordinals}

In this section we introduce notion of \textit{provable ordinals} for subsystems of $\sfZ_2$, and define $\ord(T_0)$ for a subsystem $T_0$ of $\sfZ_2$. These form the basis of Gentzen-style proof theory (see \cite[Section IX.5]{Simpson2009}), which has been used to obtain various independence results in subsystems of $\sfZ_2$ (see \cite[Remark IX.5.11]{Simpson2009}). In particular, it is used in showing that Hilbert's Basis Theorem is not provable in subsystems $T_0$ of $\sfZ_2$ with $\ord(T_0)=\omega^\omega$, as in Theorem \ref{thm: RCA does not prove hilbert basis theorem}.

 \begin{definition} Let $T_0$ be a subsystem of $\sfZ_2$ which includes $\RCA_0$. A \textit{provable ordinal} of $T_0$ is a countable ordinal $\alpha$ such that, for some primitive recursive well-ordering $W \subseteq \bN$, $\abs{W} = \alpha$ and $T_0$ proves $\WO(W)$, i.e. that, via a pairing function (see \cite[\S II.2]{Simpson2009}), $W$ codes a well-ordering of length $\alpha$. The supremum of the provable ordinals of $T_0$ is denoted $\ord(T_0)$.
\end{definition}
 We now state the provable ordinals of the systems we are concerned with.
\begin{theorem}[See Theorems IX.5.4 and IX.5.7 of \cite{Simpson2009}]\label{thm: provable ordinals} We have
 \[ \ord(\RCA_0) = \ord(\WKL_0) = \omega^\omega, \]
 and
 \[ \ord(\ACA_0) = \varepsilon_0 := \sup(\omega,\omega^\omega, \omega^{\omega^\omega},\ldots). \]
 \end{theorem}

\section{Reverse Mathematical analysis}\label{Sec: Rev Math Analysis}
We first discuss results involving algebraic closures, as an excellent example highlighting the relative strengths of $\RCA_0$, $\WKL_0$, and $\ACA_0$ -- see Theorems \ref{thm: ctble field alg closure}, \ref{thm: uniqueness of ctble field alg closure}, and \ref{thm: strong algebraic closure}. We then analyse Hilbert's Nullstellensatz (Section \ref{sec: Nullstellensatz}) and Basis Theorem (Section \ref{sec: Basis Theorem}) from a reverse mathematical perspective, drawing from results in \cite{FSS1983}, \cite{SakamotoTanaka2004}, \cite{Simpson1988}, and \cite{Simpson2015}. 

\subsection{Motivating Example: Algebraic Closures}\label{sec: algebraic closure}
Unless specified otherwise, the following definitions are made in $\RCA_0$ and the results are stated over $\RCA_0$ as well.

We first define the notion of a countable field in $\RCA_0$ (refer \cite[Section II.9]{Simpson2009}. The more general definition of a countable structure in $\RCA_0$ can be found in \cite[Section 2]{FSS1983}.

\begin{definition}\label{def: ctble field}
    A \textit{countable field} $F$ consists of a set $\abs{F} \subseteq \bN$, together with binary operations $+_F,\cdot_F$, a unary operation $-_F$, and distinguished elements $0_F,1_F$ such that the system $\innerprod{\abs{F},+_F,-_F,\cdot_F,0_F,1_F}$ obeys the usual field axioms.
\end{definition}

\begin{definition}[See \S2 of \cite{FSS1983}] $\RCA_0$ proves that for any countable field $F$ and any $m \in \bN$, there exists a countable commutative ring $F[x_1,\ldots,x_m]$ consisting of $0$ along with all (G\"odel numbers of) expressions of the form 
\[ f(x_1,\ldots,x_m) = \sum_{i_1+\cdots i_m\leq n} a_{i_1\ldots i_m} x_1^{i_1}\cdots x_m^{i_m},\]
where $(i_1,\ldots,i_m) \in \bN^m$, $m \in \bN$, $a_{i_1\ldots i_m} \in K$, and $a_{i_1\ldots i_m} \neq 0$ for at least one $(i_1,\ldots,i_m) \in \bN^m$ with $i_1+\cdots+i_m = n$. This is the \textit{ring of polynomials in $m$ commuting indeterminates $x_1,\ldots,x_m$ over $F$.}   
\end{definition}

\begin{definition}
    A countable field $F$ is said to be \textit{algebraically closed} if for all nonconstant polynomials $f(x) \in F[x]$, there exists $a \in F$ such that $f(a)=0$.\footnote{By $f(a)$, we mean taking the alternative notation $\sum_{i=0}^n a_ix^i$ for the polynomial $f(x)$, and plugging in $a$ in place of $x$.}
\end{definition}

\begin{definition}\label{def: algebraic closure} Let $F$ be a countable field. An \textit{algebraic closure} of $F$ consists of an algebraically closed countable field $K$, together with a monomorphism $h: F \to K$ such that for all $b\in K$, there exists a nonconstant polynomial $f(x) \in F[x]$ such that $h(f)(b) = 0$.
\end{definition}

The following results emphasize the relative strengths of $\RCA_0$, $\WKL_0$, and $\ACA_0$, in this context of algebraic closures of countable fields. More specifically, $\RCA_0$ only proves the existence of algebraic closures; $\WKL_0$ goes one step further, proving the uniqueness of algebraic closures; and $\ACA_0$ proves the existence of \textit{strong algebraic closures}.

\begin{lemma}[Lemma II.9.3 of \cite{Simpson2009}]\label{lem: rca_0 proves QE for acf} The following are provable in $\RCA_0$:
\begin{enumerate}
    \item $\mathsf{ACF}$, i.e. the first-order theory of algebraically closed fields, admits quantifier elimination: For any formula $\varphi$ there exists a quantifier-free formula $\psi$, containing no new free variables (see \ref{sec: arithmetical hierarchy}), such that $\ACF \proves (\varphi \liff \psi)$.\footnote{The symbol $\proves$ is read as ``proves;" and we write $T\proves \varphi$ if there is a \textit{proof} of $\varphi$ from $T$. One may refer to \cite[\S 2]{Marker2010} for a rigorous definition of \textit{proof}. However, for this discussion, it suffices to have an intuitive understanding of the notion of a \textit{proof}.}
    \item For any quantifier-free formula $\varphi$, if $\mathsf{ACF} \proves \varphi$, then $\mathsf{AF} \proves \varphi$, where $\mathsf{AF}$ is the theory of fields.
\end{enumerate}
\end{lemma}

Simpson notes that these well-known results in Lemma \ref{lem: rca_0 proves QE for acf} have purely syntactical proofs\footnote{Simpson's proof relies on Tarski's syntactical quantifier elimination methods -- presented in, for example, \cite{KK1967}.} which can be transcribed in $\RCA_0$. Friedman, Simpson, and Smith then use Lemma \ref{lem: rca_0 proves QE for acf} to prove the following result.

\begin{theorem}[Theorem 2.5 of \cite{FSS1983}]\label{thm: ctble field alg closure} $\RCA_0$ proves that every countable field has an algebraic closure.\footnote{A different proof of this theorem can also be found in \cite[Theorem II.9.4]{Simpson2009}.}
\end{theorem}

Theorem \ref{thm: ctble field alg closure} is a powerful result that provides the machinery needed to prove some formulations of Hilbert's Nullstellensatz in $\RCA_0$, which we discuss in detail in Section \ref{sec: Nullstellensatz}. We now discuss two strengthenings of Theorem \ref{thm: ctble field alg closure}, which are no longer provable in $\RCA_0$.

The following theorem states that, over $\RCA_0$, WKL is equivalent to the existence of unique algebraic closures. Since WKL is itself not provable in $\RCA_0$, it follows that $\RCA_0$ does not prove the uniqueness of algebraic closures.

\begin{theorem}[Theorem 3.3 of \cite{FSS1983}]\label{thm: uniqueness of ctble field alg closure}
The following assertions are equivalent over $\RCA_0$:
\begin{enumerate}
    \item Weak K\"onig's Lemma
    \item Every countable field has a unique (up to isomorphism) algebraic closure.
\end{enumerate}
\end{theorem}

Simpson defines the notion of a strong algebraic closure as follows (see \cite[Section III.3]{Simpson2009}).

\begin{definition}
    Let $F$ be a countable field. A \textit{strong algebraic closure} of $F$ is an algebraic closure $h : F \to K$ (see Definition \ref{def: algebraic closure}) with the further property that $h$ is an isomorphism of $F$ onto a subfield of $K$. 
\end{definition}

\begin{theorem}[Theorem III.3.2 of \cite{Simpson2009}]\label{thm: strong algebraic closure} The following assertions are pairwise equivalent over $\RCA_0$:
\begin{enumerate}
    \item $\ACA_0$
    \item Every countable field has a strong algebraic closure.
    \item Every countable field is isomorphic to a subfield of a countable algebraically closed field.
\end{enumerate}
\end{theorem}

\subsection{Hilbert's Nullstellensatz}\label{sec: Nullstellensatz}

As stated before, there have been numerous formulations of the Nullstellensatz (see footnote \ref{footnote: formulations of Nullstellensatz}), many of which are hard to compare in subsystems of $\sfZ_2$. In this section we discuss some formulations, as coded in the literature, in the reverse mathematical context. 

Friedman, Simpson, and Smith describe two formulations of the Nullstellensatz in \cite[Section 2]{FSS1983}, which we label $\mathsf{HN}_1$ and $\mathsf{HN}_2$: Let $F$ be a countable field and $f_1,\ldots,f_m \in F[x_1,\ldots,x_n]$.
    \begin{itemize}
        \item[$\mathsf{HN}_1$.] $f_1,\ldots,f_m$ have a common root in some extension of $F$ if and only if $f_1,\ldots,f_m$ have a common root in an algebraic extension of $F$.
        \item[$\mathsf{HN}_2$.] $f_1,\ldots,f_m$ have a common root in some extension of $F$ if and only if $1 \notin (f_1,\ldots,f_m)$.\footnote{Note that $\mathsf{HN}_2$ is very closely related to what was stated as the Weak Nullstellensatz in Corollary \ref{Thm: Nullstellensatz weak og}.}
    \end{itemize}

\begin{theorem}[Section 2 of \cite{FSS1983}]\label{thm: RCA0 proves HN1} $\RCA_0$ proves $\mathsf{HN}_1$.
\end{theorem}
This follows from Lemma \ref{lem: rca_0 proves QE for acf}, i.e. $\RCA_0$ proving quantifier-elimination for $\ACF$, and Theorem \ref{thm: ctble field alg closure}, i.e. $\RCA_0$ proving the existence of algebraic closures for countable fields.

\begin{theorem}[Section 2 of \cite{FSS1983}]\label{thm: RCA0 proves HN2}
    $\RCA_0$ proves $\mathsf{HN}_2$.
\end{theorem}
Although this version is customarily reduced to $\HN_1$ by extending the ideal generated by $f_1,\ldots,f_m$ to a prime or maximal ideal, the general results on the existence of prime or maximal ideals require more than $\RCA_0$ (see Theorems \ref{thm: WKL_0 equivalences} and \ref{thm: ACA_0 equivalences}). The proof of $\mathsf{HN}_2$ in $\RCA_0$ relies on a method of elimination, satisfied by Kronecker's Elimination, which is explained in more detail in \cite{FSS1983} and \cite{Stevens2023}.

\begin{corollary}[Section 2 of \cite{FSS1983}]\label{cor: cor of HN2}$ \RCA_0$ proves that $g$ vanishes on $\mathbb{V}(f_1,\ldots,f_m)$ if and only if there exists $r\geq 1$ such that $g^r \in \innerprod{f_1,\ldots,f_m}$.\footnote{Note that Corollary \ref{cor: cor of HN2} is the special case of what was stated as Hilbert's Nullstellensatz in Theorem \ref{Thm: Nullstellensatz og}, when $I = \innerprod{f_1,\ldots,f_m}$, i.e. the ideal generated by $f_1,\ldots,f_m$, whose existence is provable in $\RCA_0$ (see \cite[Lemma 2.10]{FSS1983}).}
\end{corollary}

It is a consequence of Corollary \ref{cor: cor of HN2}, by elimination, that $\RCA_0$ proves the existence of the set of all such $g$'s, i.e.\ the radical ideal\footnote{More generally, given an ideal $I$ of a ring $R$, the radical ideal generated by $I$ is $\rad(I):= \set{a \in R \mid \exists n (a^n \in I)}$.\label{footnote: radical ideal}} generated by $f_1,\ldots,f_m$.

Sakamoto and Tanaka introduce what they call Hilbert's Nullstellensatz for complex numbers in \cite[Section 1]{SakamotoTanaka2004}, which we label $\mathsf{HN}_3$:

\begin{itemize}
    \item[$\mathsf{HN}_3$] For any $n,m \in \bN$, if $p_1,\ldots,p_m \in \bC[x_1,\ldots,x_n]$ have no common zeros, then $1 \in \innerprod{p_1,\ldots,p_m}$.\footnote{Note that $\HN_3$ is the special case of the backwards direction of the Weak Nullstellensatz in Corollary \ref{Thm: Nullstellensatz weak og}, when $I=\innerprod{p_1,\ldots,p_m}$.}
\end{itemize}

\begin{theorem}[Theorem 8 of \cite{SakamotoTanaka2004}]\label{thm: RCA0 proves HN3} $\mathsf{HN}_3$ is provable in $\RCA_0$.
\end{theorem}


\subsection{Hilbert's Basis Theorem}\label{sec: Basis Theorem}

Simpson codes the Basis Theorem in terms of what he calls Hilbertian rings, defined in this section. In standard set theory, like $\ZFC$, it is provable that a ring $R$ is Hilbertian if and only if every ideal in $R$ is finitely generated, i.e. $R$ is Noetherian (see Definition \ref{def: Noetherian Ring}). On the other hand, in the reverse-mathematical context, although this equivalence holds over stronger base theories like $\ACA_0$ \cite[Remark 2.2]{Simpson1988}, it does not necessarily hold over weaker theories. In particular, over $\RCA_0$, the notion of Hilbertian is somewhat stronger than every ideal being finitely generated (see \cite[Remark 2.2]{Simpson1988}). Thus, over $\RCA_0$, the Basis Theorem coded in terms of Hilbertian rings implies the version coded in terms of every ideal being finitely generated.

\begin{definition}\label{def: Hilbertian} Within $\RCA_0$, let $R$ be a countable commutative ring. We say that $R$ is \textit{Hilbertian} if for every sequence $(r_i)_{i\in \bN}$ of elements of $R$, there exists $k \in \bN$ such that for all $j \in \bN$, there exist $s_0,\ldots,s_k \in R$ such that $r_j = \sum_{i \leq k} s_i r_i$. 
\end{definition}

Simpson then gives two, provably equivalent over $\RCA_0$, formulations of the Basis Theorem (see \cite{Simpson1988}), which we label $\mathsf{HB}_1$ and $\mathsf{HB}_2$.  

\begin{itemize}
    \item[$\mathsf{HB}_1$.] For all $m \in \bN$ and all countable fields $K$, the commutative ring\\ $K[x_1,\ldots,x_m]$ is Hilbertian.
    \item[$\mathsf{HB}_2$.] For each $m \in \bN$, there exists a countable field $K$ such that the commutative ring $K[x_1,\ldots,x_m]$ is Hilbertian.
\end{itemize}

For the remainder of this section, we list results that narrow where the Basis Theorem fits in the Friedman-Simpson Hierarchy.

Recall from Theorem \ref{thm: provable ordinals} that $\ord(\RCA_0) = \omega^\omega$. That is, $\RCA_0$ can prove that any ordinal $\alpha < \omega^\omega$ is well-ordered, and cannot do so for $\omega^\omega$. The following result asserts that $\RCA_0$ does not prove $\WO(\omega^\omega)$.

\begin{lemma}[Proposition 2.6(2) of \cite{Simpson1988}]\label{lem: RCA0 not prove WO omega^omega}
    $\RCA_0 \nproves \WO(\omega^\omega)$. That is, there is no primitive recursive well-ordering $W$ of order-type $\omega^\omega$, such that $\RCA_0 \proves \WO(W)$.
\end{lemma}

\begin{lemma}[Theorem 2.7 of \cite{Simpson1988}]\label{lem: HB equiv to WO omega^omega} The following assertions are pairwise equivalent over $\RCA_0$
\begin{itemize}
    \item[1.] $\mathsf{HB}_1$
    \item[2.] $\mathsf{HB}_2$
    \item[3.] $\WO(\omega^\omega)$
\end{itemize}
\end{lemma}

In this regard, $\omega^\omega$ is a measure of what Simpson calls the ``intrinsic logical strength" of the Basis Theorem (see \cite[Section 1]{Simpson1988}).\footnote{In \cite{KY2016}, Kreuzer and Yokoyama show that many principles of first-order arithmetic, previously only known to lie strictly between $\Sigma^0_1$-Induction and $\Sigma^0_2$-Induction, are equivalent to $\WO(\omega^\omega$). Hence, these have the same ``intrinsic logical strength" as the Basis Theorem. They argue that, in some sense, $\WO(\omega^\omega)$ should be considered a \textit{natural} first-order principle between $\Sigma^0_1$-Induction and $\Sigma^0_2$-Induction, and should have its own place (see \cite[Figure 3]{KY2016}) in the Paris-Kirby Hierarchy (as defined in \cite{HP1998}).\label{footnote: Principles in between ISigma1 and ISigma2}}

Since $\mathsf{HB}_1$ and $\mathsf{HB}_2$ are provably equivalent over $\RCA_0$, we will henceforth refer to them as $\mathsf{HB}$.

\begin{theorem} \label{thm: RCA does not prove hilbert basis theorem}
    $\RCA_0$ does not prove $\mathsf{HB}$.
\end{theorem}

This is a direct consequence of the two preceding lemmas, i.e. Lemma \ref{lem: RCA0 not prove WO omega^omega} and Lemma \ref{lem: HB equiv to WO omega^omega}: Since $\RCA_0 \nproves \WO(\omega^\omega)$ and $\RCA_0 \proves (\WO(\omega^\omega) \liff \mathsf{HB})$, we may conclude that $\RCA_0 \nproves \mathsf{HB}$.

Furthermore, Simpson showed in \cite{Simpson2015} that $\RCA_0 + \Sigma^0_2$-Induction was enough to prove $\WO(\omega^\omega)$:

\begin{lemma}[Theorem 2.2 of \cite{Simpson2015}]\label{lem: WO(omega^omega) provable in RCA0 + ISigma2} $\WO(\omega^\omega)$ is provable in $\RCA_0 + \Sigma^0_2$-Induction.
\end{lemma}

Recall that $\Sigma^0_2$-Induction (see footnote \ref{footnote: Sigma^0_2-induction}) is not provable in $\RCA_0$ (see \cite[Chapter II]{Simpson2009}). Thus, $\RCA_0 + \Sigma^0_2$-Induction is strictly stronger in terms of provability than $\RCA_0$.

\begin{theorem}\label{thm: HB provable in RCA0 + ISigma2}
    $\HB$ is provable in $\RCA_0 + \Sigma^0_2$-Induction.
\end{theorem}

This is a direct consequence of Lemmas \ref{lem: HB equiv to WO omega^omega} and \ref{lem: WO(omega^omega) provable in RCA0 + ISigma2}: Since $\RCA_0 \proves (\WO(\omega^\omega) \liff \mathsf{HB})$ and $\RCA_0 + \Sigma^0_2$-Induction $\proves \WO(\omega^\omega)$, we may conclude that $\RCA_0 + \Sigma^0_2$-Induction $\proves \HB$.

Simpson also shows that the provability result in Lemma \ref{lem: WO(omega^omega) provable in RCA0 + ISigma2} is not reversible, i.e. $\RCA_0 + \WO(\omega^\omega)$ does not prove $\Sigma^0_2$-Induction \cite[Section 4]{Simpson2015}. In fact, even assuming additional machinery, namely the $\Sigma^0_2$-Bounding Principle (see \cite[Definition 2.1]{Simpson2015}), does not suffice \cite[Corollary 4.3]{Simpson2015}. This reinforces Kreuzner and Yokoyama's argument (see footnote \ref{footnote: Principles in between ISigma1 and ISigma2}) that $\WO(\omega^\omega)$ should be considered a natural first-order principle between $\Sigma^0_1$-Induction and $\Sigma^0_2$-Induction.

It is worth making explicit that for a fixed $m \in \bN$, $\RCA_0$ proves $\WO(\omega^m)$ (see \cite[Proposition 2.6(1)]{Simpson1988}). As a consequence, $\RCA_0$ proves the following
\begin{itemize}
\item[$\mathsf{HB}(m)$.] For all countable fields $K$, the commutative ring $K[x_1,\ldots,x_m]$ is Hilbertian.
\end{itemize}

In particular, $\mathsf{HB}(1)$ is what was stated as Theorem \ref{thm: Basis Theorem og}. Hence, if, by ``Hilbert's Basis Theorem" one means $\mathsf{HB}(1)$, then Hilbert's Basis Theorem is, indeed, provable in $\RCA_0$.

Moreover, the nonconstructive content in $\HB$ truly comes from what was referred to as a ``straightforward inductive argument" in Section \ref{Sec: Bg Comm Alg}. More explicitly:

\begin{theorem}[See Proposition 2.6(3) of \cite{Simpson1988}] $\RCA_0$ proves that \[\HB\ \liff\ \forall m\ \HB(m).\]
\end{theorem}

In conclusion, although some formulations of the Nullstellensatz are provable in $\RCA_0$ (see Theorems \ref{thm: RCA0 proves HN1}, \ref{thm: RCA0 proves HN2}, \ref{thm: RCA0 proves HN3}), it follows from Simpson's results (Theorems \ref{thm: RCA does not prove hilbert basis theorem} and \ref{thm: HB provable in RCA0 + ISigma2}) that the Basis theorem, for arbitrary $m \in \bN$, needs strictly more machinery. In this sense, the Basis theorem is strictly more nonconstructive than the Nullstellensatz.

\section{Additional Remarks}\label{sec: real closures}\label{sec: additional remarks}

One could study analogous results to those in Section \ref{Sec: Rev Math Analysis}, in the context of \textit{real closed fields} and \textit{formally real fields} (see \cite[Definition II.9.5 and Remark II.9.8]{Simpson2009}).\footnote{These structures are studied reverse-mathematically in \cite{FSS1983} and \cite{Simpson2009}, and model-theoretically in \cite{JL1989} and \cite{Marker2010}.} In particular, along the lines of Section \ref{sec: algebraic closure}, $\RCA_0$ proves the existence and uniqueness (up to unique isomorphism) of real closures for countable ordered fields\footnote{It is worth noting that $\RCA_0$ was not strong enough to prove the uniqueness of algebraic closures (see Section \ref{sec: algebraic closure}).} \cite[Theorems 2.12 and 2.18]{FSS1983}; $\WKL_0$ proves that every countable formally real field is orderable, and has a real closure \cite[Theorem 3.5]{FSS1983}; and $\ACA_0$ proves that every formally real field has a strong real closure \cite[Theorem III.3.2]{Simpson2009}.

In the spirit of real closed analogs, one may wonder if there is a real version of the Nullstellensatz. Unfortunately, $\bR$, the field of real numbers, not being algebraically closed poses a nontrivial challenge. Recall that for an algebraically closed field $K$, Hilbert's Nullstellensatz establishes a one-to-one (order-reversing) correspondence between the posets of radical ideals (see footnote \ref{footnote: radical ideal}) in $K[x_1,\ldots,x_n]$ and affine varieties (Definition \ref{def: variety (common zero set)} should suffice for this discussion) in $K^n$. This correspondence fails over $\bR$. For instance, the ideal $I=\innerprod{x^2+1}$ is a radical ideal in $\bR[x]$, since $\bR[x]/I \cong \bC$ is a field. The real variety $\mathbb{V}(I) = \es$ coincides with the variety $\mathbb{V}(\innerprod{1})$ determined by the ideal $\innerprod{1}$ in $\bR[x]$. Thus, two different radical ideals define the same variety in $\bR^n$, and Hilbert's Nullstellensatz fails.

There is, however, a somewhat analogous statement in the case of real closed fields (see \cite[\S 3.4]{Marker2010}).

\begin{itemize}
    \item[$\mathsf{RN}$] Let $F$ be a real closed field, and $I$ be a prime ideal in $F[x_1,\ldots,x_n]$. Then $\mathbb{V}_F(I)$ is nonempty if and only if whenever $p_1,\ldots,p_m \in F[x_1,\ldots,x_n]$ and $\sum p_i^2 \in I$, then for each $i$, $p_i \in I$. 
\end{itemize}

A proof of this Real Nullstellensatz (over $\ZFC$) can be found in \cite{Dickmann1985}. Sakamoto and Tanaka state that it would be an interesting question to decide whether the Real Nullstellensatz is provable in $\RCA_0$; however, it is unclear how this result would even be coded in $\RCA_0$ \cite[\S 4]{SakamotoTanaka2004}.

\bibliography{References}
\bibliographystyle{alpha}
\end{document}